%% file: root.tex
\documentclass[conference]{ieeeconf}  % Comment this line out if you need a4paper

\input{packages}
\input{math_commands}

\IEEEoverridecommandlockouts                              % This command is only needed if 
\title{\LARGE \bf
Maximum Entropy Differential Dynamic Programming
}

\author{Oswin So, Ziyi Wang and Evangelos A. Theodorou
\thanks{The authors are with the Autonomous Control and Decision Systems Laboratory, Georgia Institute of Technology, Atlanta, GA, USA. Email correspondance to: \href{mailto:oswinso@gatech.edu}{\texttt{oswinso@gatech.edu}}}}

\begin{document}

\maketitle
\thispagestyle{empty}
\pagestyle{empty}

%%%%%%%%%%%%%%%%%%%%%%%%%%%%%%%%%%%%%%%%%%%%%%%%%%%%%%%%%%%%%%%%%%%%%%%%%%%%%%%%
\begin{abstract} 
In this paper, we present a novel maximum entropy formulation of the Differential Dynamic Programming algorithm and derive two variants using unimodal and multimodal value functions parameterizations. 
By combining the maximum entropy Bellman equations with
a particular approximation of the cost function, we are able to obtain a new
formulation of Differential Dynamic Programming which is able to escape from
local minima via exploration with a multimodal policy. To demonstrate the efficacy of the proposed algorithm, we provide experimental results using four systems on tasks that are represented by cost functions with multiple local minima and compare them against vanilla Differential Dynamic Programming. Furthermore, we discuss connections with previous work on the linearly solvable stochastic control framework and its extensions in relation to compositionality.
\href{https://youtu.be/NHr9Kj_jnAI?utm_source=arxiv&utm_medium=arxiv&utm_id=arxiv}{\color{blue}{Link to Video}}.
\end{abstract}

\IEEEpeerreviewmaketitle

%%%%%%%%%%%%%%%%%%%%%%%%%%%%%%%%%%%%%%%%%%%%%%%%%%%%%%%%%%%%%%%%%%%%%%%%%%%%%%%%
% 6 pages for the paper + all figures, any number of pages for references, no appendix.
% =============================================================================================
\section{Introduction}
\input{sections/intro}

% =============================================================================================
\section{Maximum Entropy Bellman Equation} \label{sec: bellman}
Standard discrete-time deterministic optimal control problems minimize the cost over time horizon $(0, 1, \cdots, T)$
\begin{equation} \label{eq:er_hjb:deterministic_objective}
J(u) \coloneqq \Phi(x_T) + \sum_{t=0}^{T-1} l_t(x_t,u_t),
\end{equation}
where $l_t$ and $\Phi$ are the running and terminal costs respectively.
The state and control trajectories, $(x_t)_{t=0,\cdots,T}, x_t \in \Rb^{n_x}$ and $ (u_t)_{t=0,\cdots,T-1}, u_t \in \Rb^{n_u}$, satisfy deterministic dynamics
\begin{equation}\label{eq:er_hjb:dynamics}
    x_{t+1} = f( x_t, u_t ).
\end{equation}

In this work, we take a \textit{relaxed control} approach and consider a stochastic control policy $\pi_t(u_t|x_t)$ with the same deterministic dynamics as in \eqref{eq:er_hjb:dynamics}.
In addition, we introduce an entropy term to the original objective \eqref{eq:er_hjb:deterministic_objective}
\begin{equation} \label{eq:er_hjb:cost_function}
    J_{\pi} \coloneqq \ExP{\pi}{  \Phi(x_T) + \sum_{t=0}^{T-1} \Big( l_t( x_t, u_t ) - \alpha H \big[\pi_t \big] \Big) },
\end{equation}
where $\alpha > 0$ is an inverse temperature term, the expectation is taken with respect to $u\sim\pi(\cdot|x)$,
and $H[\pi]$ is the Shannon entropy of $\pi$ defined as
\begin{equation}
H[\pi] = -\ExP{\pi}{\ln \pi} = -\int \pi(u) \ln \pi(u) \du.
\end{equation}
For this problem formulation and the standard value function definition of $V(x)=\min_\pi J(x, \pi)$, the Bellman equation takes the form
\begin{equation} \label{eq:er_hjb:raw_bellman}
V(x) = \inf_\pi \Big\{ \Eb_\pi \big[ l(x, u) +  V'( f(x,u) ) \big] - \alpha H\big[ \pi \big] \Big\},
\end{equation}
In \eqref{eq:er_hjb:raw_bellman} and below we omit the time index $t$ for nonterminal times for simplicity and use $V'(f(x,u))$ to denote the value function at the next timestep.

Solving \eqref{eq:er_hjb:raw_bellman} results in a Gibbs distribution for $\pi^*$
\cite{wang2020variational, kim2020hamilton}.
The form of $\pi^*$ and $V$ are presented in the following lemma.
\begin{restatable}[]{lemma}{optpolval} \label{lemma:er_hjb:opt_pol_val}
The optimal policy $\pi^*$ solving the infimum in \eqref{eq:er_hjb:raw_bellman} is the Gibbs distribution
\begin{equation} \label{eq:er_hjb:optimal_policy}
    \pi^*( u | x ) = Z^{-1} \exp\left( -\frac{1}{\alpha} \Big[ V'(f(x,u)) + l(x, u) \Big] \right),
\end{equation}
where $Z$ denotes the partition function
\begin{equation} \label{eq:er_hjb:partition_fn}
    Z(x) = \int \exp\left( -\frac{1}{\alpha} \Big[ V'(f(x,u)) + l(x, u) \Big] \right) \du.
\end{equation}
Consequently, the value function $V$ takes the form
\begin{equation} \label{eq:er_hjb:value_partition_function}
    V(x) = -\alpha \ln Z(x).
\end{equation}
\end{restatable}
We refer the readers to \Cref{app:proof:opt_pol_val} in \cite{so2021maximum} for a proof of \cref{lemma:er_hjb:opt_pol_val}.

% =============================================================================================
\section{Maximum Entropy DDP} \label{sec: me_ddp}
We will now use \ac{DDP} to solve the \ac{MEOC} problem and derive the \ac{ME-DDP} algorithm.
For notational simplicity, we will drop the second-order approximation of the dynamics as in \ac{iLQR} in our description of \ac{DDP}.
The dropped second-order dynamics terms can easily be added
back in the derivations below.
We refer readers to \cite{mayne1966second, li2004iterative} for a detailed overview of the vanilla \ac{DDP} and \ac{iLQR} algorithms.

The \ac{DDP} algorithm consists of a forward pass and a backward pass.
The forward pass simulates the dynamics forward in time obtaining a set of nominal state and control trajectories $(\bar{x}_{0:T}, \bar{u}_{0:T-1})$,
while the backward pass solves the Bellman equation with a 2nd order approximation of the costs and dynamics
equations around the nominal trajectories.
The boundary conditions for the value function $V$ are obtained by performing a 2nd order Taylor expansion of
the terminal cost $\Phi$:
\begin{equation} \label{eq:er_ddp:terminal_quadratic}
\begin{aligned}
    V_{xx,T} &= \Phi_{xx}, & V_{x,T} &= \Phi_{x}, & V_T= &= \Phi.
\end{aligned}
\end{equation}
To derive the backward pass, we first perform a quadratic approximation of the cost function around $(\bar{x}, \bar{u})$
\begin{equation*} \label{eq:er_ddp:cost_approx}
    l(x, u) \approx l(\bar{x}, \bar{u})
    + \begin{bmatrix}l_x \\ l_u\end{bmatrix}\T
      \begin{bmatrix}\delta x \\ \delta u\end{bmatrix} \\
    + \frac{1}{2}
      \begin{bmatrix}\delta x \\ \delta u\end{bmatrix}\T
      \begin{bmatrix}l_{xx} & l_{xu} \\ l_{ux} & l_{uu} \end{bmatrix}
      \begin{bmatrix}\delta x \\ \delta u\end{bmatrix},
\end{equation*}
where $\delta x \coloneqq x - \bar{x}$, $\delta u \coloneqq u - \bar{u}$.
We also perform a linear approximation of the dynamics $f$:
\begin{equation*} \label{eq:er_ddp:dynamics_approx}
    f(x, u) \approx f(\bar{x}, \bar{u}) + f_x\T \delta x + f_u\T \delta u .
\end{equation*}
Define $Q \coloneqq V'(f(x,u)) + l(x,u)$, with subscripts denoting partial derivatives.
The next lemma describes the optimal policy and value function.
%%%%%%%%%%%%%%%%%%%%%%%%%%%%%%%%%%%%%%%%%%%
\begin{restatable}[]{lemma}{meddppolval} \label{lemma:me_ddp:pol_val}
The optimal policy for the approximated problem is Gaussian with mean $\delta u^*$
and covariance $\alpha Q_{uu}^{-1}$
\begin{equation}
    \pi^*(\delta u | \delta x) = \mathcal{N}(u; \delta u^*, \alpha Q_{uu}^{-1} )
\end{equation}
where the mean $\delta u^*$ has the same form as in vanilla DDP
\begin{equation} \label{eq:er_ddp:delta_u_star}
    \delta u^* = -Q_{uu}^{-1} \Big( Q_{ux} \delta x + Q_u \Big) = K \delta x + k.
\end{equation}
Consequently, the value function has the form
\begin{equation} \label{eq:me_ddp:value_fn}
    V(x)
    = \bar{V}(x) + V_{H}(\bar{x}) + V_x(\bar{x})\T \delta x + \frac{1}{2} \delta x\T V_{xx}(\bar{x}) \delta x,
\end{equation}
where
\begin{align}
    \bar{V}(\bar{x}) &= \bar{V}'(\bar{x}) + l(\bar{x}, \bar{u})
    - \frac{1}{2} Q_u\T Q_{uu} Q_u, \\
    V_{H}(\bar{x}) &= \frac{1}{2}  \Big(
        \ln \abs{ Q_{uu} } - n_u \ln( 2 \pi \alpha )
    \Big), \\
    V_x(\bar{x}) &= Q_x + K\T Q_{uu} k + K\T Q_u + Q_{ux}\T k, \\
    V_{xx}(\bar{x}) &= Q_{xx} + K\T Q_{uu} K + K\T Q_{ux} + Q_{ux}\T K.
\end{align}
\end{restatable}
%%%%%%%%%%%%%%%%%%%%%%%%%%%%%%%%%%%%%%%%%%%
\begin{remark}
Note that the update rules for $\bar{V}, V_x$ and $V_{xx}$ are exactly the same as in
vanilla DDP. 
The only difference here is the addition of the $V_{H}$ term due to the additional entropy regularization. Consequently, this term can be ignored if the $V(x)$ is not needed.
This term disappears as $\alpha \to 0$ when the problem reverts to the vanilla case.
\end{remark}
%%%%%%%%%%%%%%%%%%%%%%%%%%%%%%%%%%%%%%%%%%%
We refer the readers to \Cref{app:proof:meddp_pol_val} in \cite{so2021maximum} for a proof of \cref{lemma:me_ddp:pol_val}.

% =============================================================================================
\section{Multimodal Maximum Entropy DDP} \label{sec: mme_ddp}
While a unimodal Gaussian policy is able to achieve better exploration compared to the deterministic policy,
it is often the case that multiple modes need to be explored simultaneously to converge
to the global minimum \cite{haarnoja2017reinforcement}.
In this section, we derive a multimodal extension to the \ac{ME-DDP} introduced.

\begin{figure}[t]
    \centering
    \includegraphics[width=0.95\columnwidth]{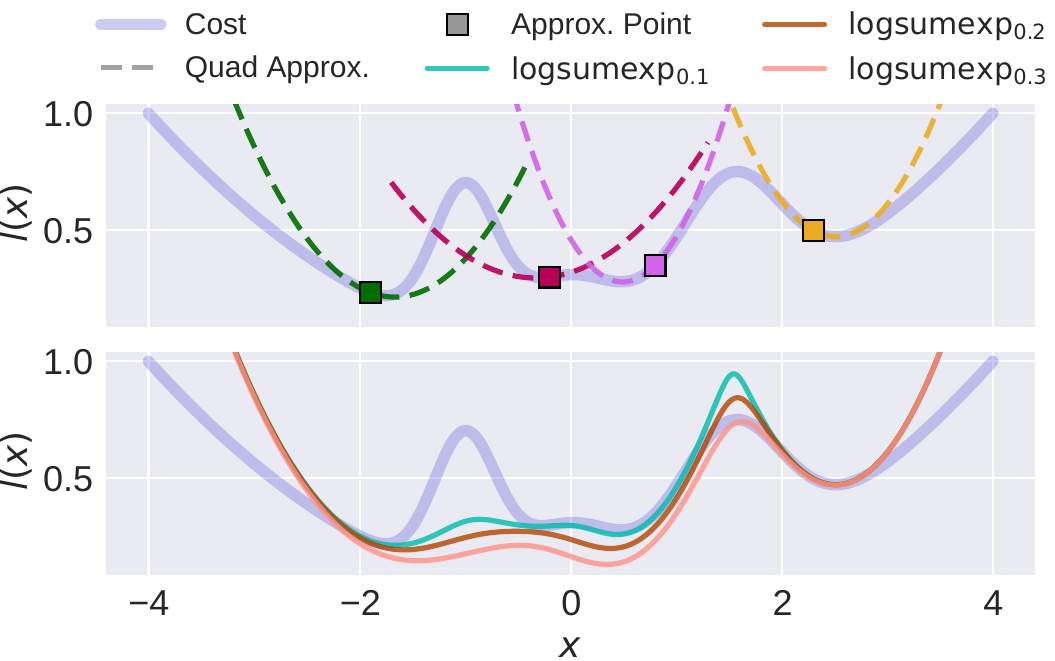}
    \caption{Comparison of the individual quadratic approximation (top) and the
    log-sum-exp approximation (bottom) of the cost function $l(x)$ with varying
    choices of inverse temperature $\alpha$. Higher $\alpha$ leads to smoother approximation.}
    \label{fig:mme_ddp:logsumexp}
\end{figure}

Let $\{ \bar{x}^{(n)}, \bar{u}^{(n)} \}_{n=1}^N$ denote $N$ different nominal state and control trajectories,
and let $\Phi^{(n)}$ corresponds to the respective quadratic
approximation of $\Phi$ around $\bar{x}^{(n)}$ and $\bar{u}^{(n)}$.
Instead of using a single quadratic approximation of $\Phi$ for the terminal cost as in \eqref{eq:er_ddp:terminal_quadratic},
we use the combined approximation $\tilde{\Phi}$
\begin{equation} \label{eq:mme:log_sum_exp_v}
V_T( x ) =
    \tilde{\Phi}( x ) \coloneqq
    -\alpha \ln \sum_{n=1}^N \exp \left(
        -\frac{1}{\alpha} \Phi^{(n)}( x )
    \right).
\end{equation}
The log-sum-exp is a smoothed combination of the local quadratic approximation and approaches
$\min_n \{ \Phi^{(n)} \}$ as $\alpha \to 0$
as shown in \cref{fig:mme_ddp:logsumexp}.

The above equation becomes easier to work with when considering the \textit{exponential transform} of 
the problem.
Define $\expalph$ to be the following function
\begin{equation}
    \expalph(y) \coloneqq \exp \Big(-\frac{1}{\alpha} y \Big).
\end{equation}
We now define the \textit{reward} $r, R_T$ and \textit{desirability} $z$ as
\begin{equation}
\begin{aligned} \label{eq:er_hjb:exp_transform}
    r_t &\coloneqq \expalph(l_t),
    &r_T &\coloneqq \expalph(\Phi),
    &z &\coloneqq \expalph(V(x)).
\end{aligned}
\end{equation}
With this transformation, note that the \textit{desirability} function $z$
is exactly the partition function $Z$ from \eqref{eq:er_hjb:value_partition_function}
and is linear in both $z'$ and $r$
(denoting $z'$ for the next timestep):
\begin{subequations}
    \begin{numcases}{}
        z(x) = Z(x) = \int z'(f(x,u)) \, r(x, u) \du, \label{eq:er_hjb:exp_bellman} \\
        z_T(x_T) = r_T(x_T).
    \end{numcases}
\end{subequations}
With this transformation, the optimal policy in \eqref{eq:er_hjb:optimal_policy} has the
following elegant form
\begin{equation} \label{eq:er_hjb:exp_optimal_policy}
    \pi( u | x ) = z(x)^{-1} z'\big( f(x,u) \big) \, r(x, u).
\end{equation}
Additionally, using the desirability function $z$ to
write \eqref{eq:mme:log_sum_exp_v}, we see that $z$ has an additive structure:
\begin{equation} \label{eq:mme:terminal_z_additive}
\begin{aligned}
    z_T(x_T)&= \sum_{n=1}^N z_T^{(n)}(x_T),
    &z^{(n)}_T(x_T) &\coloneqq r^{(n)}(x_T).
\end{aligned}
\end{equation}
The following lemma now shows that this structure holds for all time.
\begin{lemma}
    Suppose that the terminal cost has the form \eqref{eq:mme:log_sum_exp_v}.
    Then, for all $t = 0, \dots, T$, the desirability function $z$ has the following \textit{additive} structure
    \begin{equation*}
        z(x) = \sum_{n=1}^N z^{(n)}(x), \quad z^{(n)}
        \coloneqq \int z'^{(n)}( f(x,u) ) \, r(x, u) \du,
    \end{equation*}
\end{lemma}
\begin{proof}
    This holds at the terminal time from \eqref{eq:mme:terminal_z_additive}.
    Suppose that $z'(x) = \sum_{n=1}^N z'^{(n)}(x)$.
    Substituting this in \eqref{eq:er_hjb:exp_bellman}, we get
\begin{equation}
\begin{split}
    z( x )
    &= \int \left( \sum_{n=1}^N z'^{(n)}(f(x,u)) \right) r( x, u ) \du, \\
    &= \sum_{n=1}^N \int z'^{(n)}( f(x,u) ) \, r( x, u ) \du, \\
    &= \sum_{n=1}^N z^{(n)}( x ).  \label{eq:er_ddp:z_additive}
\end{split}
\end{equation}
By induction, this holds for all time.
\end{proof}
\begin{remark}
Note that the form of $z^{(n)}$ is identical to the Bellman equation \eqref{eq:er_hjb:exp_bellman}.
In other words, each $z^{(n)}$ and $V^{(n)}$ is
computed in exactly the same way as $z$ and $V$
in the unimodal case but with a different terminal condition $\Phi^{(n)}$.
\end{remark}
Substituting $V$ back for $z$ in \eqref{eq:er_ddp:z_additive} yields
\begin{equation} \label{eq:mme_ddp:value_fn}
    V(x) = -\alpha \ln \sum_{n=1}^N \exp \left( -\frac{1}{\alpha} V^{(n)}(x) \right).
\end{equation}
This result makes sense intuitively\textemdash the combined value function should be related to the minimum of
the individual approximated value functions resulting from the different nominal states.

% Note that since each $V^{(n)}$ is computed in the same way as $V$, 
% by performing the same linear and quadratic approximation of the dynamics and costs
% respectively as done in \eqref{eq:er_ddp:dynamics_approx} and \eqref{eq:er_ddp:cost_approx}
% but around $( \bar{x}^{(n)}, \bar{u}^{(n)} )$ instead,
% each $V^{(n)}$ will similarly be quadratic in the states as before.

Using the linearity of the desirability function \eqref{eq:er_ddp:z_additive},
the optimal policy \eqref{eq:er_hjb:exp_optimal_policy} has the form
% To solve for the optimal policy, we can plug in the value function
% in \eqref{eq:mme_ddp:value_fn} to \eqref{eq:er_hjb:exp_optimal_policy} to obtain
\begin{align} \label{eq:mme_ddp:optimal_policy}
    \pi( u | x )
    &= z(x)^{-1} \left( \sum_{n=1}^N z'^{(n)}( f(x,u) ) \right) r(x, u), \nonumber \\
    &= \sum_{n=1}^N
        \frac{ z^{(n)}( x ) }{ z( x ) } z^{(n)}(x)^{-1} z'^{(n)}( f(x,u) ) \, r(x, u) 
    , \nonumber \\
    &= \sum_{n=1}^N w^{(n)}( x ) \pi^{(n)}( u | x ),
\end{align}
where
\begin{align*}
    \pi^{(n)}(u | x) &\coloneqq z^{(n)}( x )^{-1} z'^{(n)}( f(x,u) ) \, r( x, u ), \\
    w^{(n)}( x ) &\coloneqq z( x )^{-1} z^{(n)}( x ), \quad \sum_{n=1}^N w^{(n)}=1.
\end{align*}
This is exactly the policy obtained in the normal case, except that we consider $z^{(n)}$ instead of $z$.
Since $V^{(n)}$ is quadratic in the state, each $\pi^{(n)}$ will be
Gaussian as before:
\begin{equation}
\begin{multlined}
    \pi^{(n)}(\delta u^{(n)} | \delta x^{(n)})
    = \mathcal{N} \Big(\delta u^{(n)}; {\delta u^{(n)}}^*, \alpha ( Q_{uu}^{(n)} )^{-1} \Big)
\end{multlined}
\end{equation}
where $\delta x^{(n)} = x - \bar{x}^{(n)}$ and $\delta u^{(n)} = u - \bar{u}^{(n)}$
are now evaluated relative to the nominal trajectories for corresponding to the $n$th approximation, with ${\delta u^{(n)}}^*$ defined analogous to \eqref{eq:er_ddp:delta_u_star} but using the approximations around
$(\bar{x}^{(n)}, \bar{u}^{(n)})$.
Importantly, from the form of \eqref{eq:mme_ddp:optimal_policy}, we see that $\pi$ is a mixture of Gaussians with component weights $w^{(n)}$ computed using the quadratic approximation of the value function \eqref{eq:me_ddp:value_fn} and are \textit{adaptive to disturbances}.
We refer the readers to \Cref{app:proof:mmeddp_details} in \cite{so2021maximum} for more details.

Since both $z$ and $\pi$ are weighted sums of $z^{(n)}$ and $\pi^{(n)}$,
% which themselves are the results of the unimodal case evaluated at different nominal trajectories $(  \bar{x}^{(n)}, \bar{u}^{(n)} )$,
computing the solution to the backward pass of \ac{MME-DDP} is equivalent to solving for \ac{ME-DDP} around the $N$ different nominal trajectories
and then \textit{composing} the value functions and policies
using \eqref{eq:mme_ddp:value_fn} and \eqref{eq:mme_ddp:optimal_policy}.

% \todo{Maybe talk about we end up approximating $l$ at $n$ different points, so there's no one
% ``consistent'' definition of $l$ if we look at the backward pass at a whole.
% Maybe talk about the alternate approach, where we use all $n$ approximations of $l$ during each step,
% but then there's added computational complexity.}

\section{Connections to Existing Works} \label{sec:connection}
\input{sections/connections}

% =============================================================================================
\input{algorithms/backwards_pass}
\input{algorithms/me_algo}
\input{algorithms/mme_algo}

% Table of comparisons.
\input{cost_table}

\section{Algorithms} \label{sec: algorithm}
There are three main algorithmic issues that need to be addressed when implementing \ac{ME-DDP} and \ac{MME-DDP}.

\textbf{Forward Pass: }
The derivation in earlier sections only describes how to perform the backward pass of \ac{DDP} to 
compute the optimal Gaussian mixture policy \eqref{eq:mme_ddp:optimal_policy},
leaving the question of how to apply the new stochastic policy for the forward pass unanswered.

For \ac{ME-DDP}, we perform multiple realizations of the stochastic
policy at each timestep.
Taking $N$ realizations for each of the $T$ timesteps will result in
polynomial growth $O(T^N)$ of required samples.
Instead, we sample the entire feedforward controls from the
stochastic policy at $t=0$ and then apply the deterministic feedback
policy for times $t=1$ to $T-1$.
To handle the added multi-modality of the optimal $\ac{MME-DDP}$ policy,
we sample from a categorical distribution to determine which of the $N$ modes will be used
for the control of a particular sample trajectory,
and then sample the feedforward controls as in the \ac{ME-DDP} case.

\textbf{Convergence: }
With a stochastic policy, the cost of the trajectory after sampling may be higher
in cost than the original trajectory or even unbounded.
To guarantee the convergence of the algorithms, we draw inspirations from \cite{dong2021replica} and apply the mean deterministic controls from the mode with the smallest cost to at least one sampled trajectory.
This guarantees that the minimum cost over all $N$ samples is monotonically decreasing, preserving the convergence properties of \ac{DDP}:
\begin{restatable}[]{lemma}{quadcvg} \label{lemma:quad_cvg}
    Each iterate of \ac{ME-DDP} and \ac{MME-DDP} results in a cost that is no worse than vanilla DDP given the current best nominal control is identical.
\end{restatable}
\noindent
We refer the readers to \Cref{app:proof:quad_cvg} in \cite{so2021maximum} for the proof.
% Let $\mathcal{B}(k)$ and $\mathcal{A}(k)$ denote the cost of the $n$th mode before and after sampling
% for the $k$th iteration respectively.
% By always keeping the controls corresponding to the mode with the smallest cost unperturbed,
% we guarantee that
% \begin{equation}
%     \min_n \mathcal{A}_n(k) \leq \min_n \mathcal{B}_n(k) \leq \min_n \mathcal{B}_n(k - 1).
% \end{equation}
% This guarantees that both $\min_n \mathcal{A}$ and $\min_n \mathcal{B}$ are monotonically decreasing functions of
% $k$, preserving the convergence properties of \ac{DDP}.

\textbf{Frequency of Sampling: }
Since the weights $w^{(n)}$ for each mode for \ac{MME-DDP} are proportional
to the value function $Z^{(n)}$,
modes which have high cost are unlikely to be resampled in the next
iteration of the forward pass
even if they will converge to a more optimal local minimum given enough \ac{DDP} iterations.
% within a basin containing a better local minima which
% it will eventually converge to after enough \ac{DDP} iterations.
To alleviate this issue, we only resample the controls for each mode after
every $m$ iterations, increasing the probability of jumping out of suboptimal
local minima.

The full \ac{ME-DDP} and \ac{MME-DDP}
algorithms are presented in \cref{alg:er_ddp} and \cref{alg:mer_ddp}, along with their backward pass in \cref{alg:backward_pass},
where variables without indices denote the entire trajectory, ($x^{(n)}$ denotes $x_{0:T}^{(n)}$).
In short, both algorithms consist of keeping the lowest cost sample and sampling the rest
from the stochastic policy $\pi$ every $m$ iterations for the forward pass, then running the backward pass for each sample.
Both passes can be executed in parallel for each sample.

\begin{figure}[t]
    \centering
    \includegraphics[width=0.45\columnwidth]{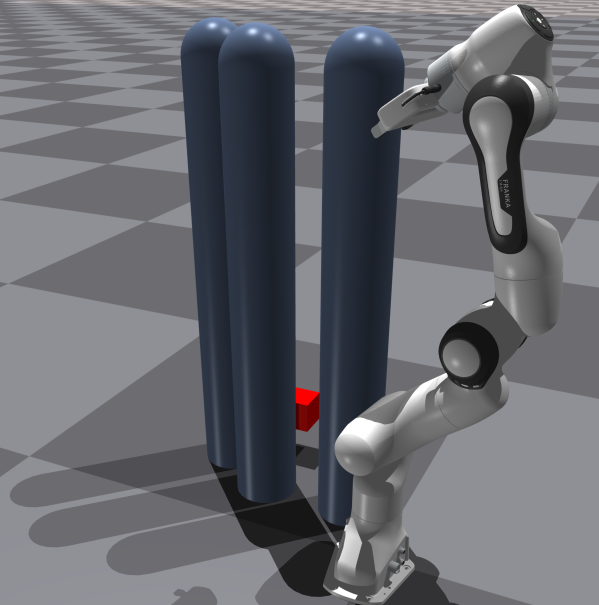}
    \caption{Task setup for the manipulator. The goal is to reach
    the red block past the obstacles while avoiding collisions.}
    \label{fig:sims:franka_render}
\end{figure}

\input{plots/state_plots}

\section{Simulations} \label{sec: results}
In this section, we compare the performance of the proposed \ac{MME-DDP} algorithm
against the \ac{ME-DDP} and vanilla \ac{DDP} algorithms on four systems:
2D Point Mass, 2D Car, Quadcopter and Manipulator.
Obstacle avoidance is implemented as a soft-constraint with
$l_{\text{obs}} = \exp \left( -\frac{ d_{ \text{obs}}^2 }{2 r_{\text{obs}}^2} \right)$,
where $d_{\text{obs}}$ and $r_{\text{obs}}$ are the distance and radius of the obstacle respectively.
The controls are zero-initialized for all systems.
For the resampling frequency, we set $m=8$.
\Cref{table:timing} compares the mean and standard deviation of the solver times
on the 2D Car problem for \ac{DDP} and \ac{MME-DDP} with $8$ modes on a Ryzen 9 3950X processor.

\input{timing_table}

\subsection{2D Point Mass}
We first test the algorithms on an illustrative 2D point-mass double integrator
reaching task while while avoiding obstacles in a maze-like environment.
Both the top and middle paths are suboptimal local minima as they are blocked
by obstacles, with the top path having an obstacle near the end of the path. The results are shown in the first row of \cref{fig:sims:2d_state_trajs}.

\subsection{2D Car}
We next test on a 2D Car with dynamics of Dubin's vehicle under jerk control.
% The system has state $x = [p_x, p_y, \theta, v, a]\T \in \Rb^5$
% and controls $u = [\omega, j]\T \in \Rb^2$,
% where $(p_x, p_y)$ describe the car's position, 
% $v, a, j$ are the linear velocity, acceleration and jerk,
% and $\theta, \omega$ are the orientation and angular velocity.
The task here is again a reaching task while avoiding two circular obstacles.
A suboptimal local minimum exists in the middle which goes in between both obstacles.
% with two more optimal minima located to the left and right of the obstacles.
The results are shown in the second row of \cref{fig:sims:2d_state_trajs}.

\subsection{Quadcopter}
We test on a 3D quadcopter with states 
$x = [p_x, p_y, p_z, \Psi, \theta, \phi, v_x, v_y, v_z, p, q, r]\T \in \Rb^{12}$
and controls $u = [f_t, \tau_x, \tau_y, \tau_z]\T \in \Rb^4$.
We refer readers to \cite{sabatino2015quadrotor} for a full description of the dynamics.
The task is to reach a target on the other side of four spherical obstacles
set up in a square pattern.
A suboptimal local minima is present in the intersection of all four
obstacles in the center which only \ac{MME-DDP} is able to consistently
escape from, as shown in the third row of \cref{fig:sims:2d_state_trajs}.

\subsection{Manipulator}
Finally, we test on a torque-controlled 7-DOF manipulator based off a simplified version of the Franka EMIKA Panda arm.
The task here is for the end effector to reach the goal position without colliding with obstacles (see \cref{fig:sims:franka_render}). 
Cylindrical obstacles are placed between the starting position and the end
effector, creating multiple suboptimal local minima in the cost landscape.
Again, only \ac{MME-DDP} is able to consistently reach the target without
intersecting any of the obstacles, as shown in the bottom row of \cref{fig:sims:2d_state_trajs}.

\textbf{Performance Comparison: }
Comparing the three algorithms,
we observe that in the case of one global minimum (Point Mass and Manipulator), vanilla \ac{DDP} gets stuck in a local minimum. Unimodal \ac{ME-DDP} explores several minima but lacks the capability to explore each sufficiently.
In contrast, the additional exploration capability helps \ac{MME-DDP} find the best minimum.
In the case of many minima of similar cost (Car and Quadcopter), vanilla and \ac{ME-DDP} get stuck in one or few suboptimal minima while \ac{MME-DDP} explores several minima simultaneously.
As obstacles are implemented as soft constraints, the shaded region around obstacles in \cref{fig:sims:2d_state_trajs} only provides a visualization and is \textbf{not} the obstacle boundary.
We also present a comparison of convergence for each algorithm in  
\cref{tab:data_table} and \cref{fig:sims:cvg_plots}.
Across all tasks, it is clear that both \ac{ME-DDP} and \ac{MME-DDP} are able
to achieve a lower mean cost than vanilla DDP due to converging to a more 
optimal local minimum.
Furthermore, \ac{MME-DDP} is able to
consistently achieve a significantly lower cost,
highlighting the advantages of multimodal exploration.

\section{Conclusion} \label{sec: conclusions}
In this paper, we derived \ac{ME-DDP} and \ac{MME-DDP}, two algorithms based off the maximum entropy formulation
of \ac{DDP} which provide improved exploration capabilities over the vanilla algorithm.
Our results suggest that the added stochasticity and multimodal exploration improves the ability
of \ac{DDP} to escape from suboptimal local minima in environments with multiple local minima.

Future work include hardware implementation to verify the exploration benefits of the proposed algorithms. On the theoretical side we will investigate the conditions and rate of convergence, as well as generalizations that include the stochastic, risk sensitive and model predictive control cases.

\FloatBarrier
%                                   % the textheight of the last page by a suitable amount.
%                                   % This command does not take effect until the next page
%                                   % so it should come on the page before the last. Make
%                                   % sure that you do not shorten the textheight too much.

%%%%%%%%%%%%%%%%%%%%%%%%%%%%%%%%%%%%%%%%%%%%%%%%%%%%%%%%%%%%%%%%%%%%%%%%%%%%%%%%

% \section*{ACKNOWLEDGMENT}

% XXX was supported by .... grant.

%%%%%%%%%%%%%%%%%%%%%%%%%%%%%%%%%%%%%%%%%%%%%%%%%%%%%%%%%%%%%%%%%%%%%%%%%%%%%%%%

% \bibliographystyle{plainnat}
\bibliographystyle{unsrtnat}
\bibliography{references}

\begin{appendices}
\include{sections/appendix}

\end{appendices}

\end{document}

%% file: packages.tex
\makeatletter
\let\NAT@parse\undefined
\makeatother

\usepackage[numbers]{natbib}
\usepackage{multicol}
\usepackage{xr-hyper}
\usepackage[bookmarks=true]{hyperref}
\usepackage{times}

\input{acronyms}
\usepackage{hyperref}
\hypersetup{colorlinks = true,
            linkcolor = blue,
            urlcolor  = black,
            citecolor = blue,
            anchorcolor = blue}

\usepackage{url}

\usepackage[utf8]{inputenc} % allow utf-8 input
\usepackage[T1]{fontenc}    % use 8-bit T1 fonts
\usepackage{url}            % simple URL typesetting
\usepackage{booktabs}       % professional-quality tables
\usepackage{amsfonts}       % blackboard math symbols
\usepackage{nicefrac}       % compact symbols for 1/2, etc.
\usepackage{microtype}      % microtypography

\usepackage{crossreftools}  % Cref in section titles.
\usepackage{thmtools}       % For restating theorems.
\usepackage{thm-restate}    % For restating theorems.

\usepackage{amsmath}
\usepackage{amssymb}
\usepackage{bm}
\usepackage{mathtools}
\usepackage{amsthm}
\usepackage{mathrsfs}
\usepackage{subcaption}
\usepackage{graphicx} % for pdf, bitmapped graphics files
\usepackage{caption}
\usepackage{cases}
\usepackage[ruled, vlined, linesnumbered, norelsize]{algorithm2e}

\usepackage[nolist]{acronym}
\usepackage{bbm}
\usepackage{xcolor}
\usepackage{wrapfig}
\usepackage{siunitx}
\usepackage{multirow}

\let\labelindent\relax
\usepackage{enumitem}

\usepackage[symbol]{footmisc}
\usepackage{float}
\usepackage[capitalise]{cleveref}
\usepackage{caption}

\usepackage{placeins}

\usepackage{tikz}
\usepackage{etoolbox}
\usepackage{siunitx}
\usepackage{booktabs}
\usepackage{dblfloatfix}

\robustify\bfseries

\usepackage{ulem}

%% file: acronyms.tex
\usepackage[nolist]{acronym}
\begin{acronym}
\acro{RL}{Reinforcement Learning}
\acro{SOC}{Stochastic Optimal Control}
\acro{MEOC}{Maximum Entropy Optimal Control}
\acro{KL}{Kullback-Leibler}
\acro{IT-MPPI}{Information Theoretic Model Predictive Path Integral}
\acro{DDP}{Differential Dynamic Programming}
\acro{iLQR}{iterative Linear Quadratic Regulator}
\acro{ME-DDP}{Maximum Entropy DDP}
\acro{MME-DDP}{Multimodal Maximum Entropy DDP}
\end{acronym}

%% file: math_commands.tex
\newcommand{\Rb}{\mathbb{R}}

\newcommand{\Eb}{\mathbb{E}}
\newcommand{\ExP}[2]{\Eb_{{#1}}{\left[#2\right]}}

\newcommand*\dif{\mathop{}\!\mathrm{d}}

\DeclareMathOperator{\expalph}{\mathcal{E}_\alpha}
\DeclareMathOperator*{\argmax}{arg\,max}
\DeclareMathOperator*{\argmin}{arg\,min}

\DeclarePairedDelimiter\abs{\lvert}{\rvert}%

\newcommand{\T}{^\mathrm{T}}

\newtheorem{lemma}{Lemma}
\newtheorem*{remark}{Remark}

\newcommand*\du{\dif u}

\SetKw{KwInParallel}{in parallel}
\newlength{\textfloatsepsave}
\newlength{\dbltextfloatsepsave}
\newlength{\floatsepsave}

%% file: sections/intro.tex
Existing methods for trajectory-optimization solve the optimization problem by iteratively relying on local information via derivatives
\cite{pontryagin1987mathematical, mukadam2016gaussian}. \ac{DDP} \cite{mayne1966second, jacobson1970differential} is a popular trajectory optimization method for nonlinear systems
used in model-based \ac{RL} and Optimal Control problems, where the problem is iteratively solved via second order
approximations of the cost and dynamics. With stage-wise positive Hessian matrices, \ac{DDP} enjoys quadratic convergence \cite{liao1991convergence}. However, these methods usually only guarantee convergence to a local minimum
and are unable to reach better local minima once converged.
In cases where there are dynamic obstacles, the cost landscape
can be highly nonconvex with  suboptimal local minima that are unsatisfactory.
Methods that try to address this problem include random restarts \cite{oleynikova2016continuous}
or via topological approaches that explicitly consider the homotopy classes of trajectories 
\cite{rosmann2017integrated, zhou2020robust}.

Maximum entropy is a technique widely used in \ac{RL} and \ac{SOC} to improve the robustness of stochastic policies. Performance robustness is achieved through an additional entropy regularization term in the cost function that improves exploration by discouraging policies from
converging to a delta distribution over the current optimal control
\cite{ziebart2010modeling, haarnoja2017reinforcement}.
In \ac{RL}, Soft Actor Critic
uses a maximum entropy objective and is considered to be state of the art for off-policy methods \cite{haarnoja2018soft}.
In \ac{SOC}, the \ac{IT-MPPI} algorithm \cite{williams2017information, wang2019information} is a generalization of this technique which uses the
forward \ac{KL} divergence between the controlled distribution and a prior distribution for regularization. The form of maximum entropy is recovered when the uniform distribution
is used as the prior distribution. In \cite{wang2020variational, lee2020generalized}, the Tsallis divergence, a generalization of the \ac{KL} divergence,
is used as a regularization term in the objective, leading to further robustness improvements.

In this paper, we consider discrete time deterministic optimal control problems and take a \textit{relaxed control} approach with entropy regularization similar to \cite{kim2020hamilton}. We propose two novel variations of \ac{DDP} under the \ac{MEOC} formulation using unimodal and multimodal Gaussian policies.
% The connection of our work to existing works optimal control compositionality law \cite{todorov2009compositionality}.
Finally, we compare the performance of both proposed algorithms against vanilla \ac{DDP} on 2D Point Mass, 2D Car, Quadcopter and Manipulator in simulation.
% \ac{iLQR} \cite{li2004iterative} is a variant of \ac{DDP} which only performs a first order approximation
% of the nonlinear dynamics and has been shown to result in equal or better convergence to \ac{DDP} in practice 
% \cite{todorov2005generalized}.

% \textbf{Contributions: }
The main contributions of this work are threefold:
\begin{itemize}
    \item We derive the Bellman equation for the discrete time
    \ac{MEOC} problem.
    \item We propose \ac{ME-DDP} and \ac{MME-DDP} to improve
    exploration over vanilla \ac{DDP}.
    \item We showcase the benefit of the improved exploration of \ac{ME-DDP}
    and \ac{MME-DDP} over vanilla \ac{DDP} in converging to better local minima
    on four different systems in simulation.
\end{itemize}

% The paper is organized as follows: In section \ref{sec: bellman}, we formulate the \ac{MEOC} problem and derive the Bellman equation. The \ac{ME-DDP} and \ac{MME-DDP} algorithms are derived in sections \ref{sec: me_ddp} and \ref{sec: mme_ddp}. The connection to compositionality theory is discussed in section \ref{sec: connection}, and a description of the proposed algorithms can be found in section \ref{sec: algorithm}. Finally, we show simulation results and conclude the paper in sections \ref{sec: results} and \ref{sec: conclusions}.

%% file: sections/connections.tex
\subsection{Compositionality and Linear Solvable Optimal Control}
A key component of our work is the compositionality of policies--solving
for the full policy $\pi$ by solving for the individual
policies $\pi^{(n)}$ then combining them via \eqref{eq:mme_ddp:optimal_policy}.
In \cite{todorov2009compositionality}, the \ac{KL} Divergence regularized control is considered
and the compositionality of controllers is introduced by exploiting the linearity of the exponential value function and the optimal policy.
Unlike \cite{todorov2009compositionality}, we allows the running cost $l(x, u)$
to be an arbitrary function of the controls $u$.
Furthermore, we provide practical algorithms in the form of \ac{ME-DDP} and \ac{MME-DDP}.
% The \ac{KL} Divergence is only used as a control cost in \cite{todorov2009compositionality}, whereas our work is
% more general and allows for the running cost $l(x, u)$ to be an arbitrary function of the
% controls $u$.
% Additionally, unlike their work, we provide practical algorithms using this compositionality
% principle in the form of \ac{ME-DDP} and \ac{MME-DDP}.
Similarly, in the field of \ac{RL} \cite{Haarnoja2018ComposableDR}, compositionality has been used on maximum entropy optimal policies to solve a conjunction of tasks by combining
maximum entropy policies which solve each of the tasks individually.

Unlike the above works, the approach our work takes focuses on the topic of exploration
rather than compositionality.
Our work is most similar to \cite{haarnoja2017reinforcement}, which shows that the multimodal
exploration
% achieved via a neural-network parametrized value function and policy network are
is able to outperform similar methods which only consider unimodal exploration policies.
However, we make use of compositionality 
and solve for the value function explicitly using \ac{DDP} which
allows for the policy to be recomputed online in realtime.

\subsection{Exploration and Control-as-Inference}
Our work is related to the Control-as-Inference framework
\cite{wang2020variational, toussaint2009robot, rawlik2012stochastic, levine2018reinforcement},
where finding the optimal policy is posed as an inference problem by minimizing the KL divergence, 
with maximum entropy emerging as a special case of this when the prior is uniform.
This framework provides a natural exploration strategy based on entropy maximization.
Since the optimal policy is usually intractable,
approaches in this area approximate the optimal policy distribution using either neural networks
or by using tractable surrogates such as Gaussian distributions. Our work can be viewed as an extension
to the latter approach by considering mixtures of Gaussians instead of unimodal Gaussians.

Our work has similar flavors to SaDDP in \cite{rajamaki2016sampled} as both methods rely on \ac{DDP} and
incorporate sampling. In their case, sampling is leveraged to address the problem of discontinuity and bypass
the use of analytical derivatives, whereas sampling is used in this work to explore multiple modes simultaneously.

%% file: algorithms/backwards_pass.tex
% \begin{algorithm}[t]
% \caption{Backwards Pass}
% \label{alg:backwards_pass}
% \begin{algorithmic}[1]
% \STATE Compute $V(T), V_x(T)$ and $V_{xx}(T)$ using $\Phi$
% \FOR{$t=T-1$ \TO $0$}
%     \STATE Compute $l$, $Q$ and their derivatives for timesteps $t$
%     \STATE Regularize $Q_{uu}$ to be PD
%     \STATE Compute $k_t, K_t, V_x(t-1), V_{xx}(t-1)$ as in Vanilla DDP
%     % \STATE $k_t \gets -Q_{uu}^{-1} Q_u$
%     % \STATE $K_t \gets -Q_{uu}^{-1} Q_{ux}$
%     % \STATE $V_x(t-1) \gets Q_x + K\T Q_{uu} k + K\T Q_u + Q_{ux}\T k$
%     % \STATE $V_{xx}(t-1) \gets Q_{xx} + K\T Q_{uu} K + K\T Q_{ux} + Q_{ux}\T K$
%     % ----
%     \STATE $\Sigma_t \gets \alpha Q_{uu}^{-1}$
%     \STATE $V_{\text{entropy}} \gets V_{\text{entropy}} + \frac{\alpha}{2} (\ln \abs{Q_{uu}} - n_u \ln(2 \pi \alpha) )$
%     % ----
% \ENDFOR
% \RETURN $k_{0:T-1}, K_{0:T-1}, \Sigma_{0:T-1}, V_{\text{entropy}, 0:T-1}$
% \end{algorithmic}
% \end{algorithm}

%
% ---- Using algorithm2e ----
%
\normalem
\begin{algorithm}[t]
\DontPrintSemicolon
\caption{Backward Pass} \label{alg:backward_pass}
Compute $V(T), V_x(T)$ and $V_{xx}(T)$ using $\Phi$\;
\For{$t = T-1$ \emph{\KwTo} $0$}
{
    Compute $l$, $Q$ and their derivatives for timesteps $t$ \;
    Regularize $Q_{uu}$ to be PD \;
    Compute $k_t, K_t, V_x(t-1), V_{xx}(t-1)$ as in Vanilla DDP \;
    $\Sigma_t \gets \alpha Q_{uu}^{-1}$ \;
    $V_{H} \gets V_{H} + \frac{\alpha}{2} (\ln \abs{Q_{uu}} - n_u \ln(2 \pi \alpha) )$ \;
}
\end{algorithm}
\ULforem

%% file: algorithms/me_algo.tex
%
% ---- Using algorithmic ----
%
% \begin{algorithm}[t]
% \caption{(Unimodal) Maximum Entropy DDP}
% \label{alg:er_ddp}
% \begin{algorithmic}[1]
% \REQUIRE Number of iterations $K$, Resample frequency $m$;
% \STATE Initialize $X^{(1:2)}, U^{(1:2)}, K^{(1:2)}, \Sigma^{(1:2)}, \pi$
% % \STATE Initialize $\pi$
% \FOR{$k=1$ \TO $K$}
%     \IF{$k \% m == 0$}
%         \STATE Set $X^{(1)}, U^{(1)}, K^{(1)}, \Sigma^{(1)}$ to the $X, U, K, \Sigma$ with the lowest cost.
%         \STATE Update $\pi$ with $X^{(1)}, U^{(1)}, K^{(1)}, \Sigma^{(1)}$
%         \STATE Sample $X^{(2)}, U^{(2)}, K^{(2)}$ from $\pi$
%     \ENDIF
%     \FOR{$n=1$ \TO $2$ \textbf{in parallel} }
%         \STATE $X^{(n)} \gets $ Rollout dynamics
%         \STATE $k^{(n)}, K^{(n)}, \Sigma^{(n)}, V_{\text{entropy}}^{(n)} \gets $ Backward Pass
%         \STATE $X^{(n)}, U^{(n)}, J^{(n)} \gets $ Line Search
%     \ENDFOR
% \ENDFOR
% % \RETURN $X^{(1:2)}, U^{(1:2)}, K^{(1:2)}$
% \end{algorithmic}
% \end{algorithm}

%
% ---- Using algorithm2e ----
%
\normalem
\begin{algorithm}[t]
\DontPrintSemicolon
\caption{(Unimodal) Maximum Entropy DDP} \label{alg:er_ddp}
\KwIn{Number of iterations $I$, Resample frequency $m$}
Initialize $x^{(1:2)}, u^{(1:2)}, K^{(1:2)}, \Sigma^{(1:2)}$ \;
\For{$i = 1$ \emph{\KwTo} $I$}
{
    \If{$i \ \% \ m = 0$}{
        $x^{(1)}, u^{(1)}, K^{(1)}, \Sigma^{(1)} \gets$ lowest code mode \;
        % Update $\pi$ with $X^{(1)}, U^{(1)}, K^{(1)}, \Sigma^{(1)}$ \;
        $x^{(2)}, u^{(2)}, K^{(2)} \sim \pi^{(1)}$ \;
    }
    \For{$n=1$ \emph{\KwTo} $2$ \emph{\KwInParallel}}
    {
        $x^{(n)} \gets $ Rollout dynamics \;
        $k^{(n)}, K^{(n)}, \Sigma^{(n)}, V_{H}^{(n)} \gets $ Backward Pass \;
        $x^{(n)}, u^{(n)}, J^{(n)} \gets $ Line Search \;
    }
}
\end{algorithm}
\ULforem

%% file: algorithms/mme_algo.tex
% %
% % ---- Using algorithmic ----
% %
% \begin{algorithm}[t]
% \caption{Multimodal Maximum Entropy DDP}
% \label{alg:mer_ddp}
% \begin{algorithmic}[1]
% \REQUIRE Number of GMM components $N$, Number of iterations $K$, Resample frequency $m$;
% \STATE Initialize $X^{(1:N)}, U^{(1:N)}, K^{(1:N)}, \Sigma^{(1:N)}, \pi$
% \FOR{$k=1$ \TO $K$}
%     \IF{$k \% m == 0$}
%         \STATE Set $X^{(1)}, U^{(1)}, K^{(1)}$ to the $X, U, K$ with the lowest cost.
%         \STATE Sample $X^{(2:N)}, U^{(2:N)}, K^{(2:N)}$ from GMM $\pi$ with weights $w^{(1:N)}$.
%     \ENDIF
%     \FOR{$n=1$ \TO $N$ \textbf{in parallel} }
%         \STATE $X^{(n)} \gets $ Rollout dynamics
%         \STATE $k^{(n)}, K^{(n)}, \Sigma^{(n)}, V_{\text{entropy}}^{(n)} \gets $ Backward Pas
%         \STATE $X^{(n)}, U^{(n)}, J^{(n)} \gets $ Line Search
%         \STATE Compute $w^{(n)}$ using $J^{(n)}$ and $V_{\text{entropy}}^{(n)}$
%     \ENDFOR
% \ENDFOR
% % \RETURN $X^{(1:N)}, U^{(1:N)}, K^{(1:N)}, w^{(1:N)}$
% \end{algorithmic}
% \end{algorithm}

%
% ---- Using algorithm2e ----
%
\normalem
\begin{algorithm}[t]
\DontPrintSemicolon
\caption{Multimodal Maximum Entropy DDP} \label{alg:mer_ddp}
\KwIn{Number of GMM components $N$, Number of iterations $I$, Resample frequency $m$}
Initialize $x^{(1:N)}, u^{(1:N)}, K^{(1:N)}, \Sigma^{(1:N)}, \pi$ \;
\For{$i = 1$ \emph{\KwTo} $I$}
{
    \If{$i \ \% \ m = 0$}{
        $x^{(1)}, u^{(1)}, K^{(1)} \gets $ lowest code mode \;
        $x^{(2:N)}, u^{(2:N)}, K^{(2:N)} \sim $ GMM $\pi$ with weights $w^{(1:N)}$ \;
    }
    \For{$n=1$ \emph{\KwTo} $N$ \emph{\KwInParallel}}
    {
        $x^{(n)} \gets $ Rollout dynamics \;
        $k^{(n)}, K^{(n)}, \Sigma^{(n)}, V_{H}^{(n)} \gets $ Backward Pass \;
        $x^{(n)}, u^{(n)}, J^{(n)} \gets $ Line Search \;
        Compute $w^{(n)}$ using $J^{(n)}$ and $V_{H}^{(n)}$ \;
    }
}
\end{algorithm}
\ULforem

%% file: cost_table.tex
% Google Spreadsheet: https://docs.google.com/spreadsheets/d/11i-FQcuoZ82O_OrsGqXX8wATfufYBOFeBMDbvQBmSN0/edit?usp=sharing
\begin{table*}[t]
    \centering
    \caption{
    Comparison of the mean and standard deviations of the cost for
    vanilla DDP, \ac{ME-DDP} and \ac{MME-DDP}, computed on 16 different DDP runs.
    The best mean cost for each system is boldfaced.
    Positive values of mean reduction correspond to a reduction. Significant reduction in mean and standard deviation can be observed from MME over both ME and vanilla DDP.
    }
    \begin{tabular}{
        @{}
        l
        S[table-format=6.2,table-number-alignment = right,table-figures-decimal=2,table-auto-round]
        S[table-format=6.2,table-number-alignment = right,table-figures-decimal=2,table-auto-round]
        S[table-format=6.2,table-number-alignment = right,table-figures-decimal=2,table-auto-round]
        S[table-format=6.2,table-number-alignment = right,table-figures-decimal=2,table-auto-round]
        S[table-format=6.2,table-number-alignment = right,table-figures-decimal=2,table-auto-round]
        S[table-format=6.2,table-number-alignment = right,table-figures-decimal=2,table-auto-round]
        S[table-format=6.2,table-number-alignment = right,table-figures-decimal=2,table-auto-round]
        S[table-format=6.2,table-number-alignment = right,table-figures-decimal=2,table-auto-round]
        @{}
    }
    \toprule
    & \multicolumn{2}{c}{Vanilla} & \multicolumn{2}{c}{ME}
    & \multicolumn{2}{c}{MME} & {MME vs Vanilla} & {MME vs ME} \\
    \cmidrule(lr){2-3} \cmidrule(lr){4-5} \cmidrule(lr){6-7} \cmidrule(lr){8-8} \cmidrule(l){9-9}
    {System} & {Mean} & {Std} & {Mean} & {Std} & {Mean} & {Std} &
    {$\Delta$Mean\%} &
    {$\Delta$Mean\%} \\
    % % End of header row.
    % ======== 2D Point Mass ========
    \midrule
    2D Point Mass
    &  32.245 & 0.0 & 10.764 & 9.550 & \bfseries 1.756 & 0.0 & 
    94.55 & 83.69 \\
    % ======== 2D Car ========
    Car
    & 5.309  & 0.0 & 4.990 & 0.643 & \bfseries 3.761 & 0.873 & 
    29.16 & 24.63 \\
    % ======== Quadcopter ========
    Quadcopter
    & 0.976  & 0.0 & 0.897 & 0.178 & \bfseries 0.536 & 0.022 & 
    45.08 & 40.25 \\
    % ======== Manipulator ========
    Manipulator
    & 22.836  & 0.0 & 20.278 & 4.679 & \bfseries 12.560 & 3.682 & 
    45.00 & 38.06 \\
    \bottomrule
    \end{tabular}
    \label{tab:data_table}
\end{table*}

%% file: plots/state_plots.tex
\begin{figure*}[t]
    \centering
    % --------------------------------------
    2D Point Mass
    \begin{subfigure}[c]{0.2\textwidth}
         \centering
         \includegraphics[trim={0 1cm 0 1cm},clip,width=0.9\textwidth]{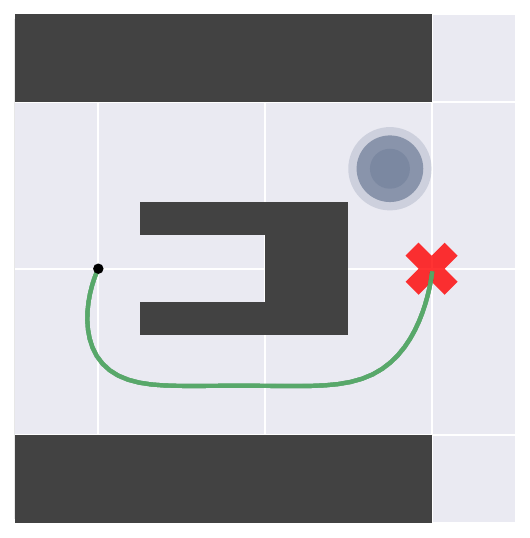}
    \end{subfigure}%
    \begin{subfigure}[c]{0.2\textwidth}
         \centering
         \includegraphics[trim={0 1cm 0 1cm},clip,width=0.9\textwidth]{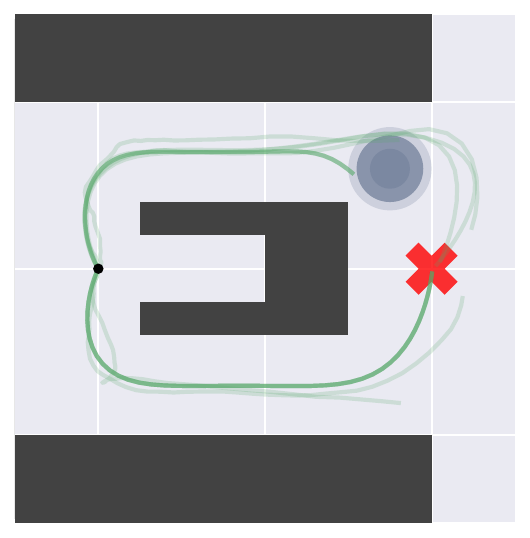}
    \end{subfigure}%
    \begin{subfigure}[c]{0.2\textwidth}
         \centering
         \includegraphics[trim={0 1cm 0 1cm},clip,width=0.9\textwidth]{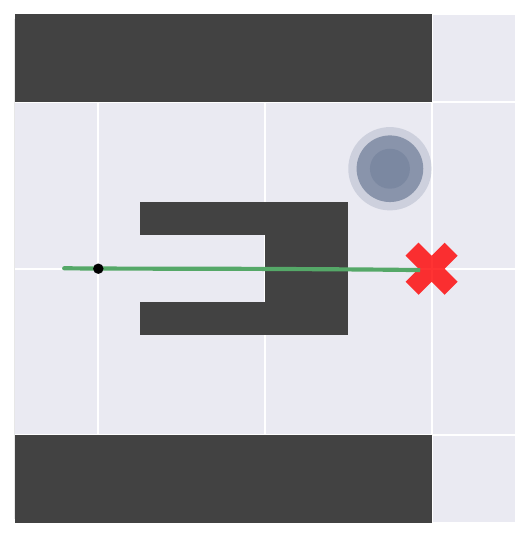}
    \end{subfigure}%
    \begin{subfigure}[c]{0.4\textwidth}
         \centering
         \includegraphics[width=0.77\textwidth]{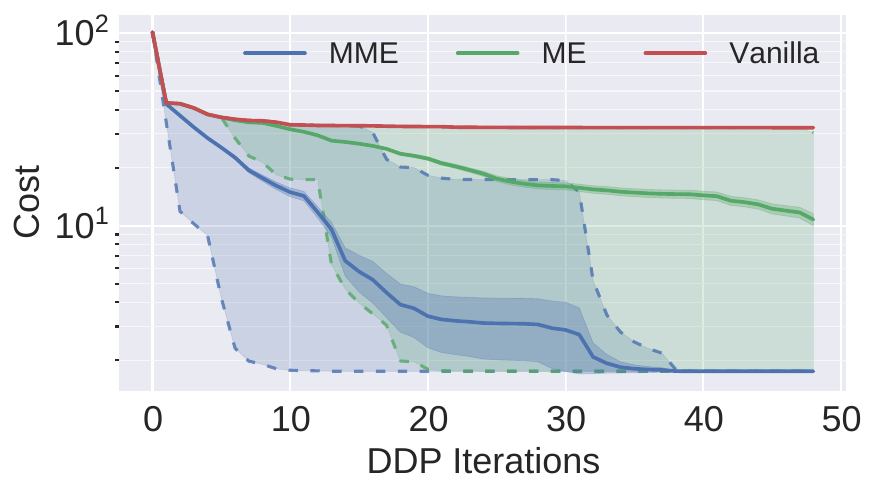}
    \end{subfigure}%
    \par\vspace{0.3em} % Vertical spacing between the rows.
    % --------------------------------------
    Car
    \begin{subfigure}[c]{0.2\textwidth}
         \centering
         \includegraphics[trim={0 1cm 0 1cm},clip,width=0.9\textwidth]{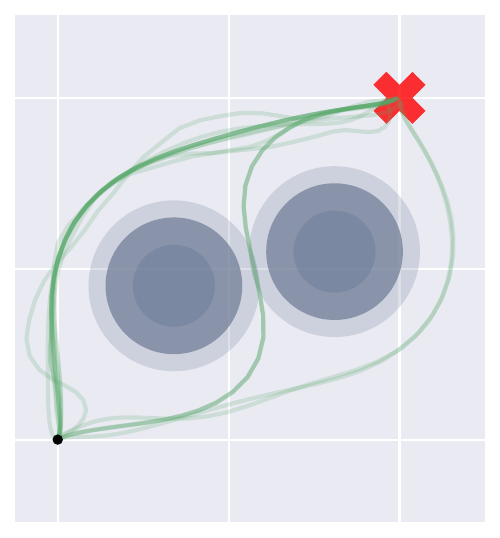}
    \end{subfigure}%
    \begin{subfigure}[c]{0.2\textwidth}
         \centering
         \includegraphics[trim={0 1cm 0 1cm},clip,width=0.9\textwidth]{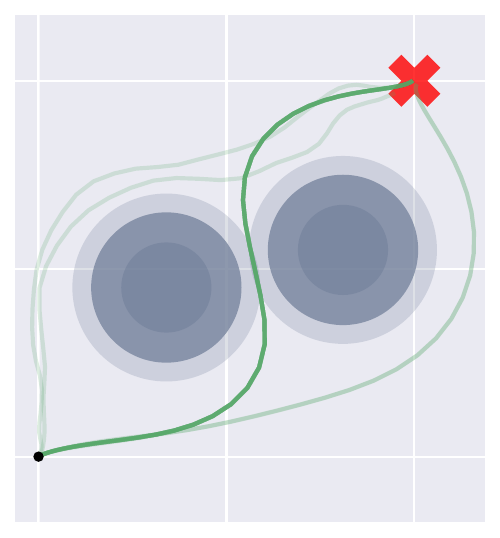}
    \end{subfigure}%
    \begin{subfigure}[c]{0.2\textwidth}
         \centering
         \includegraphics[trim={0 1cm 0 1cm},clip,width=0.9\textwidth]{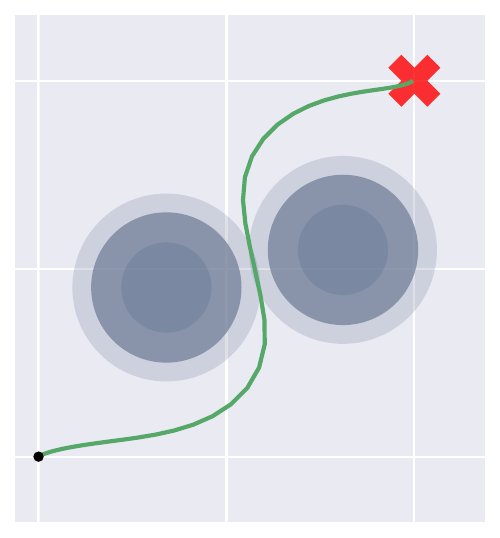}
    \end{subfigure}%
    \begin{subfigure}[c]{0.4\textwidth}
         \centering
         \includegraphics[width=0.77\textwidth]{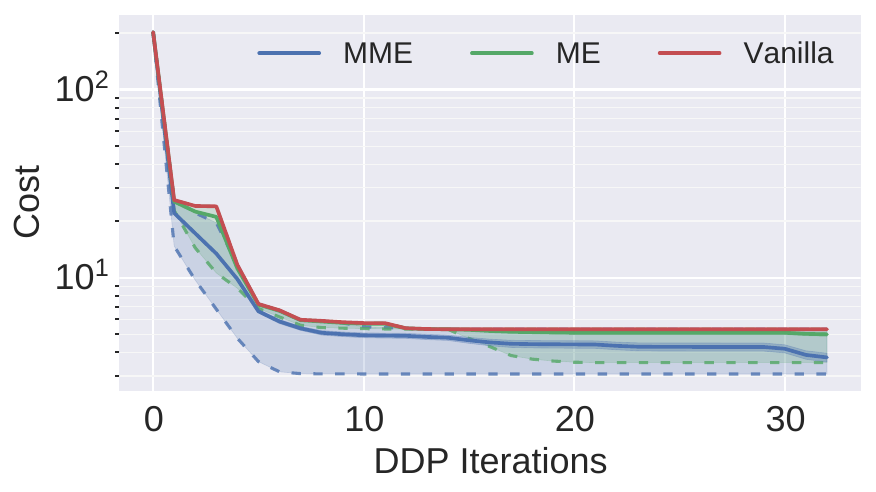}
    \end{subfigure}%
    \par\vspace{0.3em} % Vertical spacing between the rows.
    % --------------------------------------
    Quadcopter
    \begin{subfigure}[c]{0.2\textwidth}
         \centering
         \includegraphics[trim={0 1cm 0 0.2cm},clip,width=0.85\textwidth]{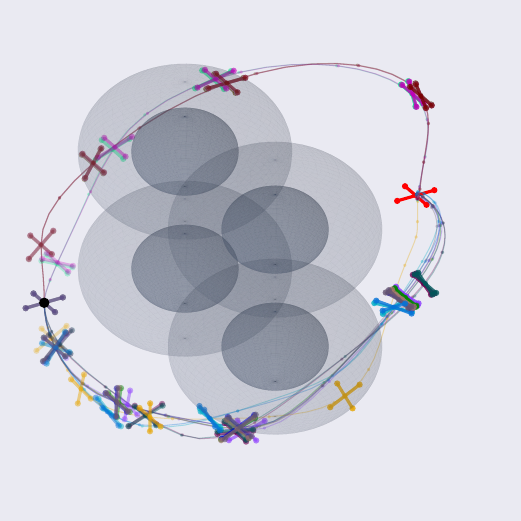}
    \end{subfigure}%
    \begin{subfigure}[c]{0.2\textwidth}
         \centering
         \includegraphics[trim={0 1cm 0 0.2cm},clip,width=0.85\textwidth]{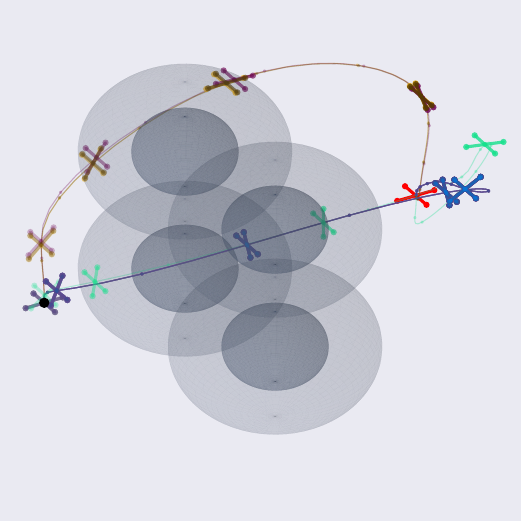}
    \end{subfigure}%
    \begin{subfigure}[c]{0.2\textwidth}
         \centering
         \includegraphics[trim={0 1cm 0 0.2cm},clip,width=0.85\textwidth]{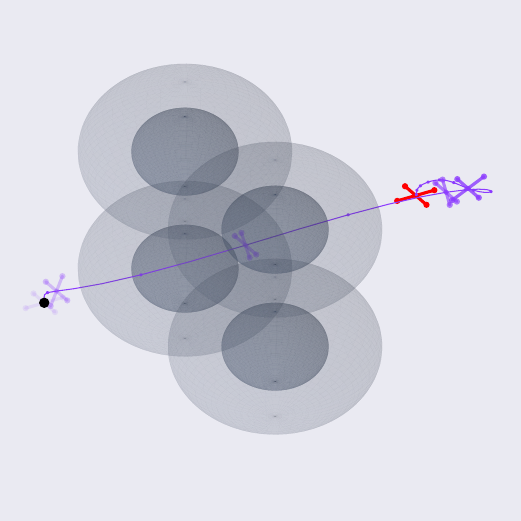}
    \end{subfigure}%
    \begin{subfigure}[c]{0.4\textwidth}
         \centering
         \includegraphics[width=0.77\textwidth]{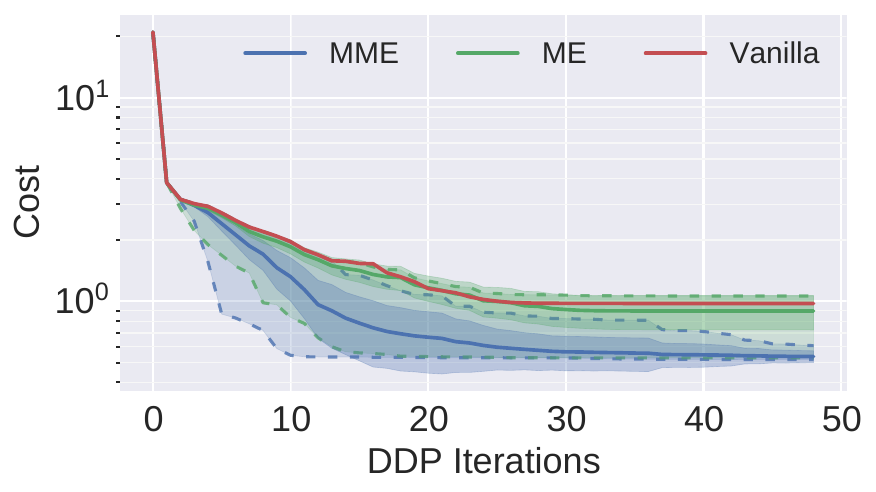}
    \end{subfigure}%
    \par\vspace{0.3em} % Vertical spacing between the rows.
    % --------------------------------------
    Manipulator
    \begin{subfigure}[t]{0.2\textwidth}
         \centering
         \includegraphics[width=0.9\textwidth]{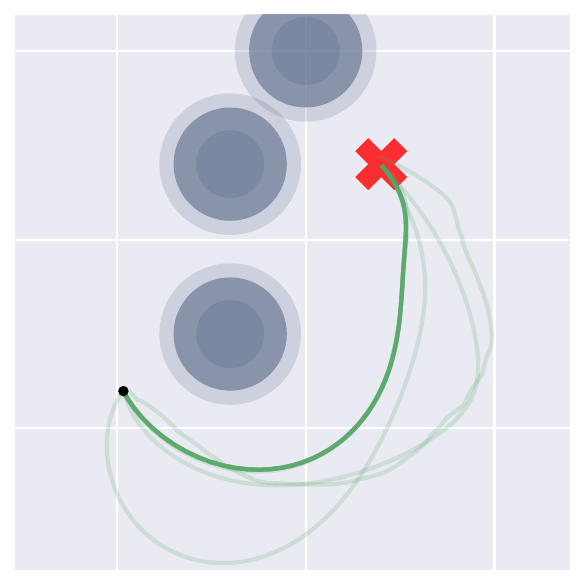}
         \caption{\ac{MME-DDP}}
    \end{subfigure}%
    \begin{subfigure}[t]{0.2\textwidth}
         \centering
         \includegraphics[width=0.9\textwidth]{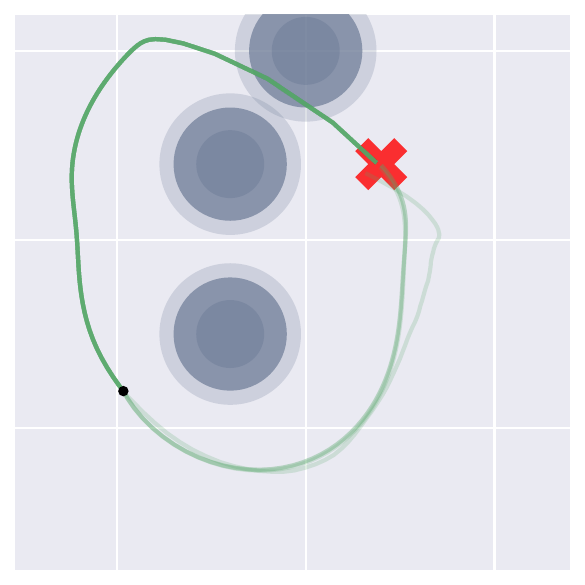}
         \caption{\ac{ME-DDP}}
    \end{subfigure}%
    \begin{subfigure}[t]{0.2\textwidth}
         \centering
         \includegraphics[width=0.9\textwidth]{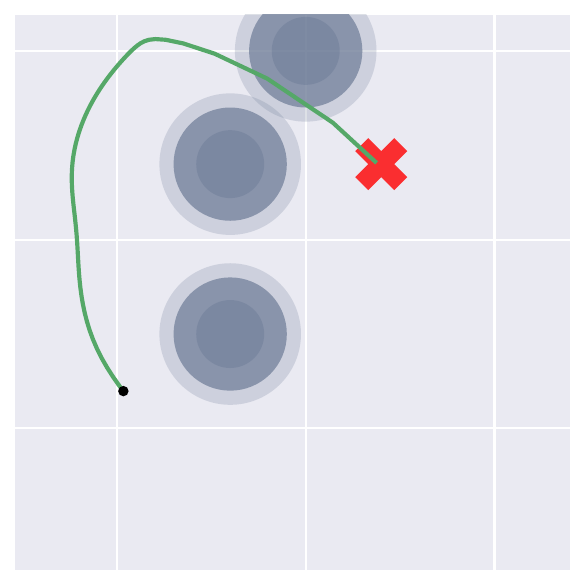}
         \caption{Vanilla DDP}
    \end{subfigure}%
    \begin{subfigure}[t]{0.4\textwidth}
         \centering
         \includegraphics[width=0.77\textwidth]{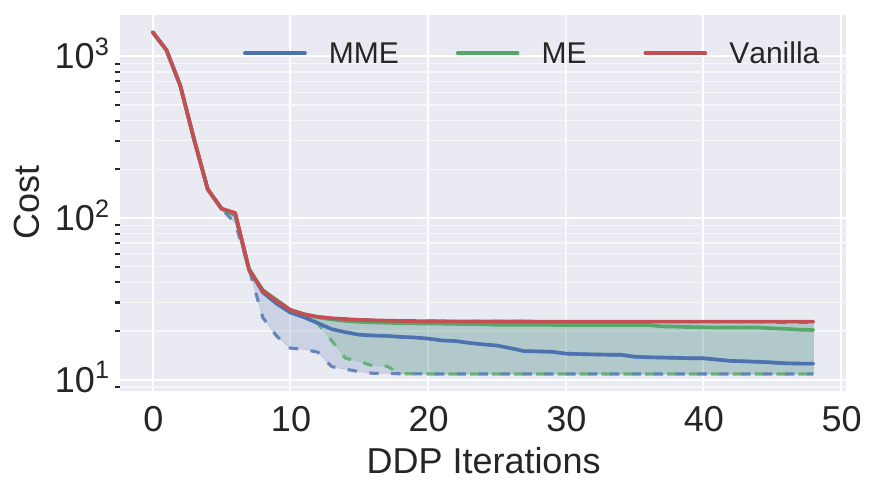}
         \caption{Convergence Plot}
         \label{fig:sims:cvg_plots}
    \end{subfigure}%
    \caption{
    (a)--(c) Position trajectories for the 2D point mass, car,
    quadcopter and manipulator systems
    from 16 different DDP runs.
    % The black and red markers represent the initial and target positions respectively
    % for each task.
    % For the quadcopter, each color represents the trajectories from a 
    % single DDP run.
    For the manipulator, the projections of the end effector trajectories on the XY-plane are plotted.
    (d) Convergence plots for all three algorithms.
    The solid line denotes the mean, the dark shaded region represents the $2\sigma$
    relative uncertainty, while the dotted lines denote the minimum and maxmimum costs.
    In all examples, \ac{MME-DDP} is able converge to a better global minimum
    due to having better exploration.
    % In all examples, without exploration, Vanilla DDP is unable to escape from the suboptimal
    % local minima.
    % \ac{ME-DDP} is not able to consistently converge to the global minima as it is only able to
    % explore a single mode.
    }
    \label{fig:sims:2d_state_trajs}
\end{figure*}

%% file: timing_table.tex
\begin{table}[t]
\centering
\begin{tabular}{@{}lrr@{}}
\toprule
Algorithm & Mean ($\si[]{\milli\second}$) & Std ($\si[]{\milli\second}$) \\ \midrule
DDP       & 1.537   & 0.016  \\ 
MMEDDP    & 2.826   & 0.118  \\ \bottomrule
\end{tabular}
\caption{Comparison of computation times for $16$ iterations averaged over 16 different runs.}
\label{table:timing}
\end{table}

%% file: sections/appendix.tex
\onecolumn

\crefalias{section}{appendix}

% Redefine section to reset equation number to 0.
\makeatletter
\let\origsection\section
\renewcommand\section{\@ifstar{\starsection}{\nostarsection}}

\newcommand\nostarsection[1]
{\sectionprelude\origsection{#1}\sectionpostlude}

\newcommand\starsection[1]
{\sectionprelude\origsection*{#1}\sectionpostlude}

\newcommand\sectionprelude{%
}

% Reset the equation number at the beginning of each section.
\newcommand\sectionpostlude{%
  \setcounter{equation}{0}
}
\makeatother

\renewcommand{\thesection}{\Alph{section}} 
\renewcommand{\thesectiondis}{\Alph{section}} 

% % \renewcommand\thesection{\Alph{section}} 
% \renewcommand\thesubsectiondis{\thesection.\Roman{subsection}}
% \renewcommand\thesubsubsectiondis{\thesubsectiondis.\alph{subsubsection}}
% \renewcommand\theparagraphdis{\arabic{paragraph}.}

% Reference equations by the Appendix number + equation number. ex. (A.1).
\renewcommand{\theequation}{\thesection.\arabic{equation}}

%%%%%%%%%%%%%%%%%%%%%%%%%%%%%%%%%%%%%%%%%%%%%%%%%%%%%%%%%%5
\section{Proof of \Cref{lemma:er_hjb:opt_pol_val}}
\label{app:proof:opt_pol_val}
\input{sections/proof_opt_pol_val}

%%%%%%%%%%%%%%%%%%%%%%%%%%%%%%%%%%%%%%%%%%%%%%%%%%%%%%%%%%5
\newpage
\section{Proof of \Cref{lemma:me_ddp:pol_val}}
\label{app:proof:meddp_pol_val}
\input{sections/proof_meddp_pol_val}

%%%%%%%%%%%%%%%%%%%%%%%%%%%%%%%%%%%%%%%%%%%%%%%%%%%%%%%%%%5
\newpage
\section{Further details on \ac{MME-DDP}}
\label{app:proof:mmeddp_details}
\input{sections/mmeddp_details}

%%%%%%%%%%%%%%%%%%%%%%%%%%%%%%%%%%%%%%%%%%%%%%%%%%%%%%%%%%5
\newpage
\section{Proof of \Cref{lemma:quad_cvg}}
\label{app:proof:quad_cvg}
\input{sections/proof_quad_cvg}

%% file: sections/proof_opt_pol_val.tex
We first restate \Cref{lemma:er_hjb:opt_pol_val} for the convenience of the reader.
\optpolval*
\begin{proof}
From the Bellman equation \eqref{eq:er_hjb:raw_bellman}, $\pi^*$ is the solution to the following
constrained optimization problem
\begin{align}
    \inf_\pi \quad &
        \ExP{u\sim\pi}{ l(x, u) + V'( f(x, u) ) } - \alpha H[ \pi(\cdot|x) ] \du \\
    \text{s.t.} \quad & \int \pi(u|x) \du = 1
\end{align}
To solve this, we first formulate the Lagrangian $\mathcal{L}$ as
\begin{align}
    \mathcal{L}
    &= \Eb_{u\sim\pi} \big[ l(x, u) + V'( f(x, u) ) \big] - \alpha H[ \pi(\cdot|x) ] 
    + \lambda \left( 1 - \int \pi(u|x) \du \right) \\
    &= \int_{\Rb^{n_u}} \pi(u|x) \Big( Q(x, u) + \alpha \ln \pi(u|x) - \lambda \Big) \du + \lambda 
\end{align}
where we have defined $Q(x,u) \coloneqq l(x,u) + V'( f(x, u) )$ for simplicity.
Applying first-order optimality conditions then gives
\begin{align}
    0 &= l(x, u) + V'( f(x, u) ) + \alpha \big( 1 + \ln \pi^*(u|x) \big) - 1 \\
    \implies \pi^*(u|x) &= \exp \left( -\frac{1}{\alpha} \Big[ l(x,u) + V'( f(x,u) ) - \lambda \Big] - 1 \right)
\end{align}
Using the constraint on the integral of $\pi^*$ then gives us
\begin{align}
1 = \int \pi^*(u | x) &= \int \exp \left( -\frac{1}{\alpha} \, Q(x,u) \right)
    \exp\left( \frac{1}{\alpha} \lambda - 1\right) \du \\
\implies
    \exp\left( \frac{1}{\alpha} \lambda - 1\right) &= \left( \int \exp \left( -\frac{1}{\alpha} \, Q(x,u) \right) \du \right)^{-1} \\
&= Z(x)^{-1}, \qquad Z(x) \coloneqq \int \exp \left( -\frac{1}{\alpha} \, Q(x,u) \right) \du
\end{align}
Hence, the optimal policy $\pi^*$ is of the form
\begin{equation} \label{eq:app:opt_pol}
    \pi^*(u|x) = Z(x)^{-1} \exp \left( -\frac{1}{\alpha} \Big[ l(x,u) + V'( f(x, u) ) \Big] \right) \du
\end{equation}
Finally, to solve for the value function $V$ in \eqref{eq:er_hjb:raw_bellman}, we plug in the optimal policy $\pi^*$ \eqref{eq:app:opt_pol} to obtain
\begin{align}
    V(x)
    &= \int \pi^*(u|x) \Big( Q(x,u) + \alpha \ln \pi^*(u|x) \Big) \du \\
    &= \int \pi^*(u|x) \Big( Q(x,u) - \alpha \ln Z(x) - Q(x, u) \Big) \du \\
    &= - \alpha \ln Z(x) \int \pi^*(u|x) \du \\
    &= - \alpha \ln Z(x)
\end{align}
\end{proof}

%% file: sections/proof_meddp_pol_val.tex
We start by restating \Cref{lemma:me_ddp:pol_val} for the benefit of the reader.
\meddppolval*
\begin{proof}
We start by performing a Taylor expansion of $Q$ with respect to $x$ and $u$ to obtain
\begin{gather}
    V'(f(x,u)) + l(x,u) \approx V'(f(\bar{x},\bar{u})) + l(\bar{x}, \bar{u}) + \delta Q, \\
    %%%%%%%%%%%%%%%%%5
    \delta Q \coloneqq
    \begin{bmatrix}Q_x \\ Q_u\end{bmatrix}\T \begin{bmatrix}\delta x \\ \delta u\end{bmatrix}
    + \frac{1}{2}
      \begin{bmatrix}\delta x \\ \delta u\end{bmatrix}\T
      \begin{bmatrix}Q_{xx} & Q_{xu} \\ Q_{ux} & Q_{uu} \end{bmatrix}
      \begin{bmatrix}\delta x \\ \delta u\end{bmatrix}. \label{eq:me_ddp:delta_Q}
\end{gather}
Expanding $\delta Q$ from \eqref{eq:me_ddp:delta_Q}, collecting terms and completing the square yields
\begin{align}
\delta Q
&= \Big( \Big[ Q_u + Q_{xu} \delta x \Big] + \frac{1}{2} \delta u\T Q_{uu} \Big)\T \delta u
    + \Big( Q_x \delta x + \frac{1}{2} \delta x\T Q_{xx} \delta x\Big) \\
&= \frac{1}{2} ( \delta u - \delta u^*)\T Q_{uu} (\delta u - \delta u^*)
    - \frac{1}{2} \delta {u^*}\T Q_{uu} \delta u^*
    + \Big( Q_x \delta x + \frac{1}{2} \delta x\T Q_{xx} \delta x\Big)
\end{align}
where $\delta u^*$ is defined as
\begin{equation}
    \delta u^* \coloneqq -Q_{uu}^{-1} \left( Q_{ux} \delta x + Q_u \right)
    = K \delta x + k.
\end{equation}
Plugging the above into \eqref{eq:er_hjb:optimal_policy} then yields
\begin{align}
    \pi^*(u|x)
    &\propto \exp\left(
        -\frac{1}{\alpha} \frac{1}{2} ( \delta u - \delta u^*)\T Q_{uu} (\delta u - \delta u^*)
    \right), \\
    &= \exp\left( -\frac{1}{2} ( \delta u - \delta u^*)\T \Sigma^{-1} (\delta u - \delta u^*) \right),
    \qquad \Sigma \coloneqq \alpha Q_{uu}^{-1},
\end{align}
which is the pdf of a multivariate Gaussian distribution with mean $\delta u^*$ and covariance $\alpha Q_{uu}^{-1}$.

Finally, to compute the value function, we first compute $Z$:
\begin{align}
    Z(x)
    &= \int \exp\left( -\frac{1}{\alpha} \left[ Q(x, u) \right] \right) \du \\
    &= 
\begin{aligned}[t]
    &\exp\left(-\frac{1}{\alpha}\left[
        V'(f(\bar{x},\bar{u})) + l(\bar{x}, \bar{u}) +
        Q_x \delta x + \frac{1}{2} \delta x\T Q_{xx} \delta x
        - \frac{1}{2}( \delta u - \delta u^*)\T \Sigma^{-1} (\delta u - \delta u^*) 
    \right]\right) \\
    &\quad \int \exp\left( -\frac{1}{2} ( \delta u - \delta u^*)\T \Sigma^{-1} (\delta u - \delta u^*) \right) \du
\end{aligned} \\
&= \exp\left(-\frac{1}{\alpha} \tilde{V}(x) \right)
    \Big( (2\pi)^{n_u} \abs{\alpha Q_{uu}^{-1} } \Big)^{\frac{1}{2}}
\end{align}
where $\tilde{V}$ contains all the terms that corresponds to the case of vanilla DDP.
Note that
\begin{align}
    -\alpha \ln \left[ \Big( (2\pi)^{n_u} \abs{\alpha Q_{uu}^{-1} } \Big)^{\frac{1}{2}} \right]
    &= -\frac{\alpha}{2} \Big( n_u \ln (2 \pi \alpha) + \ln \abs{ Q_{uu}^{-1} } \Big) \\
    &= \frac{\alpha}{2} \Big( \ln \abs{Q_{uu}} - n_u \ln (2 \pi \alpha) \Big) \\
    &\coloneqq V_H(\bar{x})
\end{align}
Consequently, applying \eqref{eq:er_hjb:value_partition_function} yields
\begin{align}
    V(x)
    &= \tilde{V}(x) - \frac{\alpha}{2} \ln \Big( (2\pi)^{n_u} \abs{\alpha Q_{uu}^{-1} }\Big) \\
    &= \bar{V}(x) + V_{H}(\bar{x}) + V_x(\bar{x})\T \delta x + \frac{1}{2} \delta x\T V_{xx}(\bar{x}) \delta x,
\end{align}
where
\begin{align}
\bar{V}(\bar{x})
    &= \bar{V}'(\bar{x}) + l(\bar{x}, \bar{u}) - \frac{1}{2} k\T Q_{uu} k, \\
    &= \bar{V}'(\bar{x}) + l(\bar{x}, \bar{u}) - \frac{1}{2} Q_u\T Q_{uu}^{-1} Q_u, \\
V_{H}(\bar{x}) &= \frac{\alpha}{2}  \Big(
        \ln \abs{ Q_{uu} } - n_u \ln( 2 \pi \alpha )
    \Big), \\
V_x(\bar{x})
    &= Q_x + K\T Q_{uu} k + K\T Q_u + Q_{ux}\T k, \\
    &= Q_x - Q_{xu} Q_{uu}^{-1} Q_u, \\
V_{xx}(\bar{x})
    &= Q_{xx} + K\T Q_{uu} K + K\T Q_{ux} + Q_{ux}\T K, \\
    &= Q_{xx} - Q_{xu} Q_{uu}^{-1} Q_{ux}.
\end{align}
\end{proof}

%% file: sections/mmeddp_details.tex
From \eqref{eq:mme_ddp:optimal_policy}, the optimal policy for \ac{ME-DDP} is a mixture distribution of Gaussians
\begin{equation} \label{eq:app:mme_opt_pol}
    \pi(u|x) = \sum_{n=1}^N w^{(n)}(x) \, \pi^{(n)}(u|x)
\end{equation}
where
\begin{align}
    \pi^{(n)}(\delta u^{(n)} | \delta x^{(n)})
    &= \mathcal{N} \Big(\delta u^{(n)}; {\delta u^{(n)}}^*, \alpha ( Q_{uu}^{(n)} )^{-1} \Big) \\
    %%%%%%%%%%%%%
    w^{(n)}( t, x_t )
    &\coloneqq z_t( t, x )^{-1} z_t^{(n)}( t, x ) \nonumber \\
    &= \frac{ \exp\left( -\frac{1}{\alpha} \left[
        V^{(n)}(t, x_t)
    \right] \right) }
    { \sum_{n'=1}^N \exp\left( -\frac{1}{\alpha} \left[ V^{(n)}(t, x_t) \right] \right) } \\
    %%%%%%%%%%%%%%
    V^{(n)}(t, x_t)
    &= \bar{V}^{(n)}(t, \bar{x}^{(n)}) 
        + V_H^{(n)}(t, \bar{x}^{(n)})
        + V_x^{(n)}(t, \bar{x}^{(n)})\T \delta x^{(n)}
        + \frac{1}{2} \delta {x^{(n)}}\T V_{xx}^{(n)}(t, \bar{x}^{(n)}) \delta x^{(n)}
\end{align}
where we have included time indices to emphasize in $w^{(n)}$ and $V^{(n)}$ to emphasize that these variables are both \textit{time} and \textit{state} dependent.
Assuming that the algorithm has converged, we should have $V_x^{(n)} \approx 0$.
Consequently, each value function $V^{(n)}$ will be centered around the nominal trajectory $\bar{x}^{(n)}$ and
serve as a ``distance metric'' weighted by the nominal cost $\bar{V}^{(n)}$ of the trajectory.
Hence, when a disturbance occurs, the weights $w^{(n)}$ of the mixture distribution $\pi$ in \eqref{eq:app:mme_opt_pol} will adapt to reflect the new best mode based on the existing quadratic approximations of each value function.

%% file: sections/proof_quad_cvg.tex
\quadcvg
\begin{proof}
Let $\bar{u}$ denote the current nominal control and let $u^{\text{DDP}}$ denote the next control iterate from vanilla DDP.

We first show that this holds for \ac{ME-DDP}.
Let $\bar{u}^{(1)}, \bar{u}^{(2)}$ denote the current nominal control trajectories for \ac{ME-DDP}, with
\begin{equation}
    \bar{u}^{(1)} = \bar{u}, \quad J( \bar{u}^{(1)} ) \leq J( \bar{u}^{(2)} )
\end{equation}
As noted in the remark for \Cref{lemma:me_ddp:pol_val}, the update rule for \ac{ME-DDP} are exactly the same as in vanilla DDP, with the mean control of $\pi^*$ matching that of vanilla DDP.
Let $\hat{u}^{(1)}$ and $\hat{u}^{(2)}$ denote the updated means after one step of \ac{ME-DDP},
such that
\begin{align}
    u^{(1)} &\gets u^{(i)}, \quad i = \argmin_{ i \in \{1, 2\} } J( \hat{u}^{(i)} ) \\
    u^{(2)} &\sim \mathcal{N}( u^{(j)}; \Sigma^{(j)} ), \quad
        j = \argmax_{ j \in \{ 1, 2 \} } J( u^{(j)} )
\end{align}
Consequently, we must have that
\begin{equation}
    \min\{ J(u^{(1)}), J(u^{(2)} ) \}
    = J( u^{(1)} )
    \leq J ( \hat{u}^{(1)} )
    = J( u^{\text{DDP}} )
\end{equation}
This same can similarly be shown for \ac{MME-DDP} assuming that $\bar{u}^{(1)} = \bar{u}$:
\begin{equation}
    \min\{ J( u^{(1)}, \dots, J( u^{(N)} ) \}
    = J( u^{(1)} )
    \leq J( \hat{u}^{(1)} )
    = J( u^{\text{DDP}} )
\end{equation}
\end{proof}